\documentstyle{amsart}
\makeatletter
\renewcommand{\subsection}{\@startsection{subsection}{2}{\z@}%
{\baselineskip}{0.5\baselineskip}{\defaultfont\bf}}
\makeatother
\def\dj{d\kern-.30em\raise1.25ex\vbox{\hrule width .3em height .03em}}
\def\Dj{D\kern-.70em\raise0.75ex\vbox{\hrule width .3em height .03em}
\kern.03em}
\newsymbol\restr 1316
\newsymbol\between 1347
\newsymbol\varnothing 203F
\newcommand{\e}{\epsilon}
\newcommand{\k}{\kappa}
\newcommand{\grten}{\mathbin{\widehat{\otimes}}}
\newcommand{\id}{\mbox{\shape{n}\selectfont id}}
\newcommand{\ad}{\mbox{\shape{n}\selectfont ad}}
\newcommand{\Sum}{\displaystyle{\sum}}
\newcommand{\im}{\mbox{im}}
\newcommand{\hor}{\mbox{\family{euf}\shape{n}\selectfont hor}}
\newcommand{\ver}{\mbox{\family{euf}\shape{n}\selectfont ver}}
\newcommand{\inv}{{i\!\hspace{0.8pt}n\!\hspace{0.6pt}v}}
\def\Mor{\mbox{Mor}}
\def\1{\varnothing}
\def\bmh{\Omega_{h\!\hspace{0.8pt}o\!\hspace{0.6pt}r}}
\def\map#1{#1_\star}
\def\bim#1{\cal{E}_{#1}}
\def\RepG{\mbox{Rep}(G)}
\def\adj{\varpi}
\newtheorem{lem}{Lemma}
\newtheorem{pro}[lem]{Proposition}
\begin{document}
\title[Quantum Principal Bundles]
{Quantum Principal Bundles as Hopf-Galois Extensions}
\author{Mi\'co \Dj ur\Dj evi\'c}
\address{Instituto de Matematicas, UNAM, Area de la Investigacion
Cientifica, Circuito Exterior, Ciudad Universitaria, M\'exico DF,
CP 04510, MEXICO}
\maketitle
\begin{abstract} It is shown that every quantum principal bundle with a
compact structure group is a Hopf-Galois extension. This property
naturally extends to the level of general differential structures, so
that every differential calculus over a quantum principal bundle with a
compact structure group is a
graded-differential variant of the Hopf-Galois extension.
\end{abstract}
\section{Introduction}
\renewcommand{\thepage}{}
The aim of this letter is to establish some connections between quantum
principal bundles and Hopf-Galois extensions. In the framework of
non-commutative differential geometry \cite{C}, the Hopf-Galois extensions
\cite{hge} are understandable as analogs of principal bundles.
If $M$ is a classical smooth manifold, $G$ a Lie group and
$P$ a classical principal $G$-bundle over $M$, then the right action of
$G$ on $P$ induces a diffeomorphism
$$ P\times G\leftrightarrow P\times_M\!P\qquad
(p,g)\mapsto(p,pg).$$
The freeness of the action of the structure group is equivalent to the
injectivity of the introduced canonical map.

Incorporating the above diffeomorphism property to non-commutative
geometry, where $M$ and $P$ are general quantum spaces (described by
noncommutative algebras) and $G$ is a
quantum group (Hopf algebra), leads to the concept of a Hopf-Galois
extension.

A more general approach to defining noncommutative-geometric concept of
a principal bundle is to generalize only the idea that $G$ acts
freely on $P$.
Such {\it quantum principal bundles} \cite{D2} trivially include
Hopf-Galois extensions, but the converse is generally not true.

In the next section we show that if $G$ is a compact matrix quantum
group \cite{W1}, then every quantum principal $G$-bundle $P$, defined as
a quantum space on which the structure group acts freely on the right,
is actually a Hopf-Galois extension.

In Section~3 the same question is
analyzed for differential calculi over quantum
principal bundles. It turns out that every differential structure
\cite{D2} over a
quantum principal bundle $P$ is understandable as a Hopf-Galois
extension, at the level of graded-differential algebras \cite{p}.

Finally, in Section~4 some concluding remarks are made.
\section{The Level of Spaces}

Let $G$ be a compact matrix quantum group, represented by a Hopf
*-algebra $\cal{A}$, the elements of which play the role of
polynomial functions on $G$.
The coproduct, counit and the antipode will be denoted by $\phi$, $\e$
and $\k$ repsectively. Let $\RepG$ be the category of finite-dimensional
unitary representations of $G$. For each $u\in\RepG$ we shall denote by
$H_u$ the carrier unitary space, so that we have the comodule structure
map $u\colon H_u\rightarrow H_u\otimes\cal{A}$. We shall denote by
$\times$ the product in $\RepG$. Let $\cal{T}$ be the
complete set of mutually non-equivalent irreducible representations of
$G$. Let $\1$ be the trivial representation of $G$, acting in
$H_\1=\Bbb{C}$.

Let $M$ be a quantum space, represented by a *-algebra $\cal{V}$. Let
$P=(\cal{B},i,F)$ be a quantum principal $G$-bundle over $M$. By
definition \cite{D2}, the map
$F\colon\cal{B}\rightarrow\cal{B}\otimes\cal{A}$ is a *-homomorphism
satisfying
$$ (\id\otimes\phi)F=(F\otimes\id)F\qquad (\id\otimes\e)F=\id.$$
It corresponds to the dualized right action of $G$ on $P$. The map
$i\colon\cal{V}\rightarrow\cal{B}$ is a *-monomorphism representing the
projection map. The image $i(\cal{V})\subseteq\cal{B}$ coincides with
the $F$-fixed point subalgebra of $\cal{B}$. Geometrically  this means
that $M$ is identificable, via the projection, with the corresponding
orbit space.
The final condition
corresponds to the requirement that $G$ acts freely on $P$. We require
that for each $a\in\cal{A}$ there exists elements $q_k,b_k\in\cal{B}$ such
that
$$\sum_k q_kF(b_k)=1\otimes a.$$
Equivalently, the map $X\colon\cal{B}\otimes\cal{B}\rightarrow
\cal{B}\otimes\cal{A}$ given by $X(q\otimes b)=qF(b)$
is surjective. This definition implies that the domain of
$X$ can be factorized to the tensor product over $\cal{V}$, which will
be denoted by $\otimes_M$. We shall
denote by the same symbol the projected map
$X\colon\cal{B}\otimes_M\!\cal{B}
\rightarrow\cal{B}\otimes\cal{A}$.

In this section we are going to prove that $X$ is bijective. In other
words, this means that every quantum principal bundle $P$ is
undersandable as a Hopf-Galois extension. The proof will be based on the
representation theory \cite{W1} of compact matrix quantum groups. We
shall construct explicitly the inverse of $X$.

For each $u\in\RepG$ let $\bim{u}$ be the space of intertwiners
between $u$ and $F$. These spaces
are $\cal{V}$-bimodules, finite and projective on both sides. We have
$\bim{\1}=\cal{V}$, in a natural manner.

The following natural decomposition holds
$$\cal{B}=\sideset{}{^\oplus}\sum_{\alpha\in\cal{T}}\cal{B}^\alpha,\qquad
\cal{B}^\alpha=\bim{\alpha}\otimes H_\alpha.$$
The above decomposition is specified by
$\varphi(x)\leftrightarrow\varphi\otimes x$.

\renewcommand{\thepage}{\arabic{page}}
As explained in \cite{D-tann}, the product map in $\cal{B}$ induces a
natural identification
$$ \bim{u}\otimes_M\!\bim{v}\leftrightarrow\bim{u\times v}$$
for each $u,v\in\RepG$. Furthermore, every intertwiner $f\in\Mor(u,v)$
induces, via the composition map, a bimodule homomorphism $\map{f}\colon
\bim{v}\rightarrow\bim{u}$, so that we have a contravariant functor from
$\RepG$ to the category of finite projective bimodules. In particular,
the contraction maps $\between^u\colon H_u^*\otimes H_u
\rightarrow\Bbb{C}$ induce bimodule injections $\map{\between^u}
\colon\cal{V}\rightarrow\bim{u}^*\otimes_M\!\bim{u}$, where we have identified
$\bim{u^{\!c}}\leftrightarrow\bim{u}^*$
with the help of the *-structure in $\cal{B}$. Here $u^c\colon
H_u^*\rightarrow H_u^*\otimes\cal{A}$ and $\bim{u}^*$ are
the corresponding conjugate representation and the conjugate module.

By construction, we have
\begin{equation}\label{Dij}
\map{\between^u}(1)=\sum_k\nu_k\otimes\mu_k\qquad
\sum_k\nu_k(e_i^*)\mu_k(e_j)=\delta_{ij}1,
\end{equation}
where $\{e_i\}$ is an orthonormal basis in $H_u$, and $\{e_i^*\}$ the
corresponding biorthogonal basis in $H_u^*$.

Let $\tau\colon\cal{A}\rightarrow\cal{B}\otimes_M\!\cal{B}$ be a linear
map defined by
\begin{equation}\label{tau}
\tau(u_{ij})=\sum_k\nu_k(e_i^*)\otimes\mu_k(e_j),
\end{equation}
where $u\in\cal{T}$. This map can be naturally extended to $\tau\colon
\cal{B}\otimes\cal{A}\rightarrow\cal{B}\otimes_M\!\cal{B}$, by imposing
the left $\cal{B}$-linearity.
\begin{pro}\label{pro1}
The maps $\tau$ and $X$ are mutually inverse.
\end{pro}
\begin{pf} A direct computation gives
$$ X\tau(b\otimes u_{ij})=b\sum_kX\bigl(\nu_k(e_i^*)\otimes\mu_k(e_j)
\bigr)=b\sum_{kn}\nu_k(e_i^*)\mu_k(e_n)\otimes u_{nj}=b\otimes u_{ij},$$
where we have used \eqref{Dij} and the intertwining property.
Similarly, if $\mu\in\bim{u}$ then
\begin{multline*}
\tau X\bigl(b\otimes\mu(e_i)\bigr)=\sum_j\tau
\bigl(b\mu(e_j)\otimes
u_{ji}\bigr)=\sum_{kj}b\mu(e_j)\nu_k(e_j^*)\otimes\mu_k(e_i)\\
=\sum_{kj}b\otimes\mu(e_j)\nu_k(e_j^*)\mu_k(e_i)=b\otimes\mu(e_i),
\end{multline*}
because of the $F$-invariance of $\Sum_j\mu(e_j)\nu_k(e_j^*)$.

In conclusion, $X$ is bijective and $X^{-1}=\tau$.
\end{pf}

It is worth noticing that if the structure group is non-compact,
there exists a variety of quantum principal bundles which are
not Hopf-Galois extensions (the map $X$ has a non-trivial kernel).
\section{The Level of Differential Structures}

In this section we shall consider general differential structures
on quantum principal bundles, and prove that they are always
understandable as Hopf-Galois extensions for differential algebras.
Considerations of this section are based on the general
theory of differential forms and connections developed in \cite{D2}.

\subsection{Preliminaries About Differential Calculi}
Let $\Gamma$ be a bicovariant *-calculus \cite{W2} over $G$.
Let $\Gamma_{\inv}$ be the space of left-invariant elements of
$\Gamma$. There exists the canonical projection map $\pi\colon
\cal{A}\rightarrow\Gamma_{\inv}$, given by the formula
$$\pi(a)=\k(a^{(1)})d(a^{(2)}). $$
The group $G$ naturally acts on the space $\Gamma_{\inv}$, via the
projection $\adj\colon\Gamma_{\inv}\rightarrow\Gamma_{\inv}\otimes
\cal{A}$ of the adjoint action. Explicitly,
$$\adj\pi=(\pi\otimes\id)\ad,\qquad\ad(a)=a^{(2)}\otimes\k(a^{(1)})a^{(3)}.$$
Let us assume that
the higher-order calculus on $G$ is described by the corresponding
\cite{D1} universal envelope $\Gamma^\wedge$.
Let $\Gamma^\otimes$ be the tensor bundle algebra over $\Gamma$. Let us
denote by $\Gamma_{\inv}^{\wedge,\otimes}$ the corresponding
left-invariant subalgebras. These algebras are *-invariant. The
action $\adj$ admits natural extensions to
$\Gamma_{\inv}^{\wedge,\otimes}$.

Let us consider a graded-differential algebra $\Omega(P)$, representing
the calculus on $P$. By definition, this means that
$\Omega^0(P)=\cal{B}$, and that $\cal{B}$ generates the whole
differential algebra $\Omega(P)$, as well as that $F$ is extendible
(necessarily uniquely) to a homomorphism $\widehat{F}\colon
\Omega(P)\rightarrow\Omega(P)\grten\Gamma^\wedge$ of graded-differential
*-algebras. Geometrically, the map $\widehat{F}$ corresponds to the
``pull back'' of the right action map.

Let $\hor(P)\subseteq\Omega(P)$ be a *-subalgebra representing
horizontal forms. By definition,
$$\hor(P)=\widehat{F}^{-1}\bigl(\Omega(P)\otimes\cal{A}\bigr).$$
We have $$F^\wedge(\hor(P))\subseteq\hor(P)\otimes\cal{A},$$
where $F^\wedge\colon\Omega(P)\rightarrow\Omega(P)\otimes\cal{A}$ is the
corresponding right action of $G$, defined as a composition
$F^\wedge=(\id\otimes p_*)\widehat{F}$, and
$p_*\colon\Gamma^\wedge\rightarrow\cal{A}$ is the projection map.

Let $\Omega(M)\subseteq\hor(P)$ be a (graded-differential)
*-subalgebra consisting of $F^\wedge$-invariant elements.
Equivalently, $\Omega(M)$ is a $\widehat{F}$-fixed point
subalgebra of $\Omega(P)$. It represents the calculus on $M$.

By definition, a connection on the bundle
$P$ is every first-order hermitian linear
map $\omega\colon\Gamma_{\inv}\rightarrow\Omega(P)$ satisfying
$$\widehat{F}\omega(\vartheta)=\sum_k\omega(\vartheta_k)\otimes c_k
+1\otimes\vartheta$$
where $\Sum_k\vartheta_k\otimes c_k=\adj(\vartheta)$. As explained in
\cite{D2}, every quantum principal bundle $P$ admits a connection.

Let us fix a splitting of the form
$\Gamma^\otimes_{\inv}=S^\wedge_{\inv}\oplus\Gamma_{\inv}^\wedge$,
realized via hermitian grade-preserving
right-covariant section
$\iota\colon\Gamma^\wedge_{\inv}\rightarrow\Gamma^\otimes_{\inv}$
of the corresponding
factor-projection map. Here $S^\wedge\subseteq\Gamma^\otimes$ is
the quadratic ideal defining $\Gamma^\wedge$, and $S_{\inv}^\wedge\subseteq
\Gamma_{\inv}^\otimes$ is its left-invariant part.

Let $\omega^\wedge\colon\Gamma^\wedge_{\inv}
\rightarrow\Omega(P)$ be the restriction of the unital multiplicative extension
$\omega^\otimes$ to $\Gamma^\wedge_{\inv}$. The following natural
decompositions hold
\begin{align}
m_\omega\colon\hor(P)\otimes\Gamma^\wedge_{\inv}&\leftrightarrow\Omega(P)
\label{dec1}\\
m_\omega^-\colon\Gamma^\wedge_{\inv}\otimes\hor(P)&\leftrightarrow
\Omega(P),\label{dec2}
\end{align}
where $m_\omega$, and $m_\omega^-$ are the unique left/right
$\hor(P)$-linear extensions of $\omega^\wedge$.

\subsection{Hopf-Galois Extensions}
The map $X$ admits a natural extension to $\widehat{X}\colon
\Omega(M)\otimes_M\!\Omega(P)\rightarrow\Omega(P)\grten\Gamma^\wedge$,
explicitly
$$\widehat{X}(\alpha\otimes\beta)=\alpha\widehat{F}(\beta).$$
In this subsection, the symbol $\otimes_M$ will be used for the tensor
product over $\Omega(M)$.
\begin{pro}\label{pro2}
The map $\widehat{X}$ is bijective.
\end{pro}
\begin{pf}
The restriction map
\begin{equation}\label{restr}
\widehat{X}\colon\Omega(P)\otimes_M\!\hor(P)\rightarrow\Omega(P)\otimes
\cal{A}
\end{equation}
is bijective, as directly follows from the decomposition
\eqref{dec2} and the proof of Proposition~\ref{pro1}.

We shall first prove that all
spaces $\Omega(P)\grten\Gamma^\wedge_k$ are included in the image of
$\widehat{X}$, using the induction on $k$. Here $\Gamma^\wedge_k=
\Sum^\oplus_*\Gamma^{\wedge j}$, where the summation is performed over
$j\leq k$. Let us fix a connection $\omega$ on $P$.

Let us assume that the statement holds for some $k$. Then for each
$\vartheta\in\Gamma^{\wedge k+1}_{\inv}$ and $\alpha\in \Omega(P)
\otimes_M\!\hor(P)$ we have
$$ \widehat{X}\bigl(\alpha\omega^\wedge(\vartheta)\bigr)
=\widehat{X}(\alpha)\vartheta+\xi,$$
where $\xi\in\Omega(P)\grten\Gamma^\wedge_k$, so that
$\widehat{X}(\alpha)\vartheta\in\im(\widehat{X})$. Hence, $\widehat{X}$ is
surjective.

We prove the injectivity of $\widehat{X}$.  Let us suppose
that $w=\Sum_k\alpha_k\omega^\wedge(\vartheta_k)$ belongs to
$\ker(\widehat{X})\setminus\{0\}$, where
$\alpha_k\in\Omega(P)\otimes_M\!
\hor(P)$ and the elements $\vartheta_k\in\Gamma^\wedge_{\inv}\setminus\{0\}$
are homogeneous and linearly independent.

This implies $\Sum^*_k\widehat{X}(\alpha_k)\vartheta_k=0$ where the
summation is performed over indexes corresponding to the elements
$\vartheta_k$ of the maximal degree. However, this is a contradiction,
and hence $\widehat{X}$ is injective.
\end{pf}

Therefore, $\Omega(P)$ is a Hopf-Galois extension of $\Omega(M)$, at the
level of graded-differential algebras \cite{p}. The map
$\widehat{X}$ intertwines the corresponding differentials,
and hence we have a cohomology isomorphism
$H\bigl(\Omega(P)\otimes_M\!\Omega(P)\bigr)
\leftrightarrow H(P)\otimes H(\Gamma^\wedge)$.

The inverse $\widehat{\tau}$ of $\widehat{X}$ can be constructed
explicitly in the following way. We have
$$\widehat{\tau}(\vartheta)=1\otimes\omega(\vartheta)-
\sum_k\omega(\vartheta_k)\Delta_k\qquad\widehat{\tau}
\restr\cal{B}\otimes\cal{A}=\tau,$$
where $\Sum_k\vartheta_k\otimes c_k=\adj(\vartheta)$ and $X(\Delta_k)=c_k$.
Using the above formulas, and the identity
$$ \widehat{\tau}(w\otimes\alpha\beta)=w\sum_{\gamma\delta}
(-1)^{\partial\alpha\partial\gamma}\gamma\widehat{\tau}(\alpha)\delta$$
where $\Sum_{\gamma\delta}\gamma\otimes\delta=
\widehat{\tau}(\beta)$, we can compute the values of $\widehat{\tau}$ on
arbitrary elements.

\section{Concluding Remarks}

The proof of Proposition~\ref{pro1} is strongly based on the
representation theory of compact matrix quantum groups. However, the
compactness assumption figures only implicitly in the proof of
Proposition~\ref{pro2}. Actually, the proof works for general quantum
structure groups. It is sufficient to ask that the bundle admits
connections, that there exists a splitting $\iota$ and that the
following natural decomposition holds
\begin{equation}\label{dech}
\Omega(M)\otimes_M\!\cal{B}\leftrightarrow\hor(P)\leftrightarrow
\cal{B}\otimes_M\!\Omega(M).
\end{equation}

In the compact case the three assumptions hold
automatically. Under the
mentioned assumptions, if $X$ is bijective then $\widehat{X}$ is
bijective, too.

The first two assumptions ensure the
existence of decompositions \eqref{dec1} and \eqref{dec2}. The third
assumption ensures that $\hor(P)$ is a Hopf-Galois extension of
$\Omega(M)$, if $\cal{B}$ is a Hopf-Galois extension of $\cal{V}$. In
other words, the restriction map \eqref{restr} will be bijective. The rest
of the proof is the same.

Proposition~\ref{pro2} holds also in the case when the higher-order
calculus on the structure group is described by the corresponding
\cite{W2} braided exterior algebra $\Gamma^\vee$. Because the general
formalism of differential forms and connections can be formulated
in the same way for braided exterior algebras.

Let $P=(\cal{B},i,F)$ be an arbitrary quantum principal $G$-bundle over
$M$. Let $\Omega(P)$ be an arbitrary differential structure describing the
calculus on $P$.
As explained in \cite{D2}, the canonical bimodule sequence
$$ 0@>>> \hor^1(P) @>>> \Omega^1(P) @>>> \ver^1(P) @>>> 0 $$
is exact. Here $\ver(P)$ is the graded-differential *-algebra describing
``verticalized'' differential forms on the bundle. We have
$\ver(P)\leftrightarrow\cal{B}\otimes\Gamma^\wedge_{\inv}$,
at the level of graded left $\cal{B}$-modules

Let $\bmh(P)\subseteq\hor(P)$ be the *-subalgebra generated by
$\cal{B}$ and $di(\cal{V})$. In the paper \cite{B1},
these forms are called horizontal.
The exactness condition introduced in \cite{B1} for
the differential calculus can be formulated as the equality
\begin{equation}\label{bm-md}
\bmh^1(P)=\hor^1(P).
\end{equation}

It is interesting to observe \cite{B2} that $P$ is a Hopf-Galois
extension if and only if \eqref{bm-md} holds for {\it universal differential
structures} (on both the bundle and the structure group). Therefore
Proposition~\ref{pro1} can be reformulated as the statement that
\eqref{bm-md} always holds for completely universal differential
calculus, if the structure group is compact.

On the other hand, for general differential structures the property
\eqref{bm-md} does not generally hold. For example, it is sufficient to
consider a non-universal calculus $\Gamma$ on $G$, and to assume that
the calculus on $P$ is universal. Another interesting phenomena is that
even if \eqref{bm-md} holds, generally $\bmh(P)\neq\hor(P)$ at the
higher-order levels. A simple example for this is given by the trivial
bundle $P=G$ over a one-point set $M$. If we assume that the calculus on
the structure group $G$ is described by the braided exterior algebra
$\Gamma^\vee$ and that $\Omega(P)=\Gamma^\wedge$, then
$$\Omega(M)=\Bbb{C}\iff \Gamma^\wedge=\Gamma^\vee.$$

\section{Acknowledgements}

It is a pleasure to acknowledge various interesting discussions with
Profs Thomasz Brzezi\'nski and Piotr Hajac, during Workshop on
Quantum and Classical Gauge Theories, Stefan Banach International Mathematical
Center, Warsaw. I would like to thank Profs Witold Kondracki and Robert
Budzy\'nski for kindly inviting me to visit Banach Center.

I am especially grateful to Prof Zbigniew Oziewicz, in particular
for inviting me to visit the Institute of Theoretical Physics in
Wroclaw, and for his kind hospitality in this beautiful town.
I would like to thank Prof Andrzej Borowiec for various interesting
discussions in Wroclaw.

\end{document}